\documentclass[11pt]{article}
\usepackage{epsfig, psfrag,graphicx, color, url}
\usepackage{geometry, amsmath, latexsym, amsfonts, amssymb}                

\newcommand{\R}{\mathbb{R}}
\newcommand{\Q}{\mathbb{Q}}
\newcommand{\Z}{\mathbb{Z}}
\newcommand{\C}{\mathbb{C}}

\newcommand{\T}{\mathbb{T}}
\newcommand{\N}{\mathbb{N}}

\newcommand{\II}{\mathfrak{I}}

\def\EEE{{\mathcal E}}

\def\GGG{{\mathcal G}}

\def\OOO{{\mathcal O}}

\def\YYY{{\mathcal Y}}

\def\UUU{{\mathcal U}}

\newtheorem{thm}{Theorem}[section]
\newtheorem{defn}[thm]{Definition}

\newtheorem{prop}[thm]{Proposition}

\newtheorem{question}[thm]{Question}

\newcommand{\qed}{\nopagebreak \begin{flushright}
      \rule{2mm}{2.5mm} \end{flushright}}

\newcommand{\id}{\mbox{\rm id}}                  
\newcommand{\cl}{\overline}                      

\newcommand{\intersect}{\cap}                    
\newcommand{\union}{\cup}                        
\newcommand{\diam}{\mbox{\rm diam}}					         

\newcommand{\interior}{\mbox{int}}     

\newcommand{\Isom}{\mbox{Isom}}                  
\newcommand{\mtwo}[4]                            
{\mbox{$\left(\begin{array}{cc}                  
#1 & #2 \\
#3 & #4 
\end{array}
\right)$}}

\newcommand{\dettwo}[4]                          
{\mbox{$\left|\begin{array}{cc}                  
#1 & #2 \\
#3 & #4 
\end{array}
\right|$}}

\newcommand{\be}{\begin{enumerate}}
\newcommand{\eb}{\end{enumerate}}
\newcommand{\bi}{\begin{itemize}}
\newcommand{\ib}{\end{itemize}}
\newcommand{\bl}{\begin{list}}
\newcommand{\lb}{\end{list}}

\newcommand{\gap}{\vspace{5pt}}                 

\newcommand{\roundness}{\mbox{\rm Round}}

\newcommand{\wtU}{{\widetilde{U}}}

\newcommand{\wtA}{\widetilde{A}}

\newcommand{\XX}{\mbox{$\mathfrak{X}$}}

\begin{document}

\title{Examples of coarse expanding conformal maps}

\author{Peter Ha\"issinsky and Kevin M. Pilgrim}

\date{\today}

\maketitle

\abstract{In previous work, a class of noninvertible topological dynamical systems $f: X \to X$ was introduced and studied; we called these {\em topologically coarse expanding conformal} systems.  To such a system is naturally associated a preferred quasisymmetry (indeed, snowflake) class of metrics in which arbitrary iterates distort roundness and ratios of diameters by controlled amounts; we called this {\em metrically coarse expanding conformal}.  In this note we extend the class of examples to several more familiar settings, give applications of our general methods, and discuss implications for the computation of conformal dimension.}

\tableofcontents

\section{Introduction} The goal of this note is threefold: first, to give further concrete examples of so-called topologically coarse expanding conformal and metrically coarse expanding conformal dynamical systems, introduced in \cite{kmp:ph:cxci}; second, to apply the general theory developed there to some recent areas of interest; and lastly, to pose some problems about the conformal gauges associated with these dynamical systems.  

\subsection*{Topologically cxc systems} Let $X$ be a compact, separable, metrizable topological space; for simplicity, we assume here that $X$ is connected
and locally connected.  
Suppose that $f: X \to X$ is a continuous, open, closed (hence surjective) map which in addition is a degree $d \geq 2$ {\em branched covering} in the sense of \cite{edmonds:fbc}. 
Let  $\UUU_0$ be a finite open cover of $X$ by open connected sets, and for $n \geq 0$  set inductively $\UUU_{n+1}$ to be the covering 
whose elements are connected components of inverse images of elements of $\UUU_n$.   

\begin{defn}[Topologically cxc]
The topological dynamical system $f: X \to X$ is {\em topologically coarse expanding conformal with respect to $\UUU_0$} provided the following axioms hold.  
\be  

\item {\bf [Exp]}  The mesh of the coverings $\UUU_n$ tends to zero as $n \to \infty$.  That is, for any finite open cover $\YYY$ of $X$ by open sets, 
there exists $N$ such that for all $n \geq N$ and all $U \in \UUU_n$, there exists $Y \in \YYY$ with $U \subset Y$.  

\item {\bf [Irred]}  The map $f: X \to X$ is {\em locally eventually onto near $X$}:  for any $x \in X$ and any neighborhood $W$ of $x$ in $X$, there is  some $n$ with $f^n(W) \supset X$.

\item {\bf [Deg]} The set of degrees of maps of the form $f^k|\wtU: \wtU \to U$, where $U \in \UUU_n$, $ \wtU \in \UUU_{n+k}$, and $n$ and $k$ are arbitrary, has a finite maximum.  
\eb
\end{defn}
It is easy to see the property of being cxc is preserved under refinement of $\UUU_0$, so this is indeed an intrinsic property of the dynamical system.
Since we have assumed $X$ to be connected, [Irred] is a consequence of [Exp].

\subsection*{Metrically cxc systems} 
Suppose now $X$ is a metric space. 
\gap 

\noindent{\bf Roundness.} Let $A$ be a bounded, proper  subset of $X$ with nonempty interior.  Given $a \in \interior(A)$, define the {\em outradius} of $A$ about $a$ as  
\[ L(A,a)=\sup\{|a-b|: b \in A\}\]
and the {\em inradius} of $A$ about $a$ as 
\[ \ell(A,a)=\sup\{ r : r \leq L(A,a) \; \mbox{ and } \; B(a,r) \subset A\}.\]
The {\em roundness of $A$ about $a$} is defined as 
\[ \roundness(A,a) = L(A,a)/\ell(A,a) \in [1, \infty).\]
\gap

\begin{defn}[Metric cxc]
The metric dynamical system $f: X \to X$ is {\em metrically 
coarse expanding conformal} provided there exist
\bi
\item continuous, increasing embeddings $\rho_{\pm}:[1,\infty) \to [1,\infty)$, the {\em forward and backward
roundness distortion functions}, and

\item increasing homeomorphisms $\delta_{\pm}:(0,\infty) \to (0,\infty)$, the {\em forward and backward relative
diameter distortion functions}
\ib
satisfying the following axioms:
\be
\setcounter{enumi}{3}

\item {\bf [Round]} for all $n\geq 0$ and $k\geq 1$, and for all
\[ U \in \UUU_n, \;\;\wtU \in \UUU_{n+k}, \;\; \tilde{y} \in \wtU, \;\;  
y
\in U\]
if
\[ f^{k}(\wtU) = U, \;\;f^{k}(\tilde{y}) = y \]
then the {\em backward roundness bound}
\[
\roundness(\wtU, \tilde{y}) <
\rho_-(\roundness(U,y))
\]
and the {\em forward roundness bound}
\[
\roundness(U, y) <
\rho_+(\roundness(\wtU,
\tilde{y}))
\]
hold.

\item {\bf [Diam]} $(\forall n_0, n_1, k)$ and for all
\[ U \in \UUU_{n_0}, \;\;U' \in \UUU_{n_1}, \;\;\wtU \in \UUU_{n_0+k},
\;\;\wtU'
\in
\UUU_{n_1+k}, \;\; \wtU' \subset \wtU, \;\; U' \subset U\]
if
\[ f^k(\wtU) = U, \;\;f^k(\wtU') = U'\]
then
\[ \frac{\diam\wtU'}{\diam\wtU} < \delta_-\left(\frac{\diam U'}{\diam
U}\right)\]
and
\[ \frac{\diam U'}{\diam U} < \delta_+\left(\frac{\diam \wtU'}{\diam
\wtU}\right).\]

\eb
\end{defn}

\subsection*{Conformal gauges} 
First, some notation.  We denote the distance between two points $a,b$ in a metric space $X$ by $|a-b|$.  Given nonnegative quantities $A, B$ we write $A \lesssim B$ if $A< C\cdot B$ for some constant $C>0$; we write $A \asymp B$ if $A \lesssim B$ and $B \lesssim A$.  

A homeomorphism $h$ between two metric spaces is {\em quasisymmetric} if there is distortion function $\eta: [0,\infty) \to [0,\infty)$ which is a homeomorphism satisfying $|h(x)-h(a)| \leq t |h(x)-h(b)| \implies |h(x)-h(a)| \leq \eta(t)|h(x)-h(b)|$.   For the kinds of spaces we shall be dealing with, this is equivalent to the condition that the roundness distortion of balls is uniform.  The {\em conformal gauge} of a metric space $X$ is the set of all metric spaces to which it is quasisymmetrically equivalent, and its {\em conformal dimension} is the infimum of the Hausdorff dimensions of metric spaces $Y$ belonging to the conformal gauge of $X$.   

The principle of the Conformal Elevator shows \cite[Thm. 2.8.2]{kmp:ph:cxci}:

\begin{thm}
\label{thm:promotion}
A topological conjugacy between metric cxc dynamical systems is quasisymmetric.
\end{thm}

Two metrics $d_1, d _2$ are {\em snowflake equivalent} if $d_2 \asymp d_1^\alpha$ for some $\alpha>0$.   We have \cite[Prop. 3.3.11]{kmp:ph:cxci}:

\begin{thm}
\label{thm:canonical_gauge}
If $f: X \to X$ is topologically cxc, then for all $\varepsilon>0$ sufficiently small, there exists a metric $d_\varepsilon$ such that:
\be
\item the elements of $\UUU_n, n\in \N$, are uniformly round;
\item their diameters satisfy $\diam\, U \asymp \exp(-\varepsilon n), U \in \UUU_n$, $n \in \N$.
\item if $f|B_\varepsilon(x, 4r)$ is injective, then $f$ is a similarity with factor $e^\varepsilon$ on $B(x,r)$.
\eb
Any two metrics satisfying (1), (2), (3) are snowflake equivalent.  The conformal gauge $\GGG$ of $(X, d_\varepsilon)$ depends only on the topological dynamical system $f: X \to X$.  
\end{thm}

It follows that the conformal dimension of $(X, d_\varepsilon)$ is an invariant of the topological conjugacy class of $f$, so it is meaningful to speak of the conformal dimension of the dynamical system determined by $f$.  

The two theorems above have extensions to cases where $X$ is disconnected; we refer the reader to \cite{kmp:ph:cxci} for details.  

Examples of metrically cxc systems include the following: 
\be
\item hyperbolic, subhyperbolic, and semihyperbolic rational maps, acting on their Julia sets equipped with the spherical metric;
\item quasiregular maps on Riemannian manifolds whose iterates are uniformly quasiregular;
\item smooth expanding maps on smooth compact manifolds, when equipped with certain distance functions.
\eb

In \S 2, we show that certain invertible iterated function systems in the plane, equipped with the Euclidean metric,  naturally yield metrically cxc systems.  Combined with Theorem \ref{thm:promotion}, 
this yields an  extension of a recent result of Ero\v{g}lu et al.

In \S 3, we consider skew products of shift maps with coverings of the circle.  
These arise naturally as subsystems of dynamical systems on the $2$-sphere.  We build, by hand, metrics for which they are cxc.   

In \S 4, we give examples of metrically cxc systems on the Sierpi\'nski carpet and on the Menger curve.  

In \S 5, we show how a so-called Latt\`es example can be perturbed to yield a continuous one-parameter family of topologically cxc maps on the sphere, and we pose some problems regarding the associated gauges.
\gap

\subsection*{Acknowledgement} The first author's work was  supported by project ANR ÒCannonÓ (ANR-06-BLAN-0366).   The second author acknowledges travel support from the NSF and from Texas A\&M university to attend the Dynamical Systems II conference, whose Proceedings comprise this volume.  We also thank Juan Rivera for useful conversations and the anonymous referee for valuable comments.

\section{Iterated function systems}

In this section, we apply our technology to the setup of Ero\v{g}lu, Rohde, and Solomyak \cite{eroglu:rohde:solomyak:qs}.  

For $\lambda \in \C$, $0 < |\lambda| < 1$ let $F_0(z)=\lambda z$, $F_1(z)=\lambda z + 1$.  The maps $F_0, F_1$ determine an {\em iterated function system} possessing a unique compact attractor
$A_\lambda = F_0(A_\lambda) \union F_1(A_\lambda)$.  The set $A_\lambda$ is invariant under the involution $s(z) = -z+(1-\lambda)^{-1}$.  Bandt \cite{bandt:nonlinearity:mandelbrot} observed that if $\lambda$ belongs to the set 
\[ \mathcal{T}=\{ \lambda : F_0(A_\lambda) \intersect F_1(A_\lambda) = \mbox{ a singleton,} \ \{o_\lambda\} \},\]   
then $F_0$ and $F_1 \circ s$ are inverse branches of a degree two branched covering $q_\lambda: A_\lambda \to A_\lambda$.  

In this case, $A_\lambda$ is also known to be a dendrite, i.e. compact, connected, and locally connected \cite{bandt:keller:simple}, and the unique branch point $o_\lambda$ of the map $q_\lambda$ is a cut-point of $A_\lambda$ and is nonrecurrent \cite[Thm. 2]{bandt:rao:topology}.   The complement of $o_\lambda$ in $A_\lambda$ is a disjoint union $A_\lambda^0, A_\lambda^1$ where $q_\lambda(o_\lambda)\in A_\lambda^1$.  
Associated to $q_\lambda$ is a combinatorial invariant, the {\em kneading sequence}, defined as the itinerary of $o_\lambda$ with respect to the partition $\{\{o_\lambda\}, A_\lambda^0, A_\lambda^1\}$ of $A_\lambda$.  

 It follows immediately that there exists a covering $\UUU_0$ by small open connected subsets of $A_\lambda$ such that $q_\lambda$ is topologically cxc with respect to $\UUU_0$.  It is easy to show that the Euclidean metric satisfies conditions (1)-(3) in Theorem \ref{thm:canonical_gauge}.  We conclude that the Euclidean metric is snowflake equivalent to the metric given by this theorem, and that the metric dynamical system given by $q_\lambda$ is metrically cxc. 

The metrics produced by Theorem \ref{thm:canonical_gauge} are obtained as 
visual distances on the boundary of a Gromov hyperbolic graph.
The discussion in the preceding paragraph shows that at least in this case, they admit more down-to-earth descriptions.  

Now suppose $f_c(z)=z^2+c$ is a quadratic polynomial for which the Julia set is a dendrite.  There is a standard way to partition the Julia set into three pieces $J_c^0, J_c^1, \{0\}$ such that $c \in J_c^1$, and one defines analogously the kneading sequence of $f_c$ to be the itinerary of the origin with respect to this partition.  
In the case when the orbit of $o_\lambda$ under iteration of $q_\lambda$ is finite, Kameyama \cite{kameyama:self-similar} showed that $q_\lambda$ is topologically conjugate to a unique quadratic polynomial $f_{c_\lambda}(z) = z^2+c_\lambda$ acting on its Julia set $J_{c_\lambda}$.    Since such polynomials are metrically cxc with respect to the Euclidean metric, Theorem \ref{thm:promotion} implies this conjugacy is quasisymmetric.  This recovers the first half  of  \cite[Theorem 1.1]{eroglu:rohde:solomyak:qs}, but does not yield the existence of an extension of this conjugacy to the Riemann sphere, which requires more work.  

More generally,  Ero\v{g}lu et. al. establish

\begin{prop} {\bf \cite[Prop. 5.2]{eroglu:rohde:solomyak:qs}}
Suppose $\lambda \in \mathcal{T}$  and $c\in\C$ is a parameter such that $J_c$ is a dendrite.  
If the kneading sequences of $q_\lambda$ and $f_c$ are identical,  
then $(A_\lambda, q_\lambda)$ and $(J_c, p_c)$ are topologically conjugate.
\end{prop}

Suppose now that $q_\lambda$ and $f_c$ satisfy the hypotheses of this proposition.  Since recurrence is a topological condition, the critical point of $f_c$ at the origin is nonrecurrent.  Since the Julia set $J_c$ is a dendrite, the map $f_c$ cannot have parabolic cycles.  By \cite[Thm. 4.2.3]{kmp:ph:cxci}, the map $f_c$ is metrically  cxc with respect to the Euclidean metric.  The map $q_\lambda$ is also metrically cxc with respect to the Euclidean metric.  Applying Theorem 1.3 to the topological conjugacy given by the preceding proposition, we obtain a stronger conclusion:

\begin{thm} Suppose $\lambda \in \mathcal{T}$  and $c\in\C$ is a parameter such that $J_c$ is a dendrite.  
If the kneading sequences of $q_\lambda$ and $f_c$ are identical,
then $(A_\lambda, q_\lambda)$ and $(J_c, p_c)$ are quasisymmetrically  conjugate.
\end{thm}

\begin{question}
Is the conformal dimension of $A_\lambda$ equal to $1$?  
\end{question}

Tyson and Wu \cite{tyson:wu:selfsimilar} show, for example, that the conformal dimension of the standard Sierpi\'nski gasket is equal to $1$ (but not realized) by exhibiting a family of explicit quasiconformal deformations through IFSs.  Can their techniques be adapted for the above IFSs?

\section{Skew products from Thurston obstructions}

In this section, we show that, associated to a certain combinatorial data, there exists a metric cxc dynamical system realizing the conformal dimension.  The type of combinatorial data arises naturally when considering topologically cxc maps $f: S^2 \to S^2$ possessing combinatorial obstructions, in the sense of Thurston, to the existence of an invariant quasiconformal (equivalently, conformal) structure; see \cite{DH1}, \cite{haissinsky:pilgrim:cxcIII}.  

Here is the outline.
\be
\item We begin with a directed multigraph $\GGG$ with weighted edges satisfying certain natural expansion and irreducibility conditions. 
\item From this data, and  a {\em snowflake parameter} $\alpha>0$, we define an associated map $g: \II_1 \to \II_0$, $\II_1 \subset \II_0$ on a family of Euclidean intervals whose inverse branches constitute a so-called {\em graph-directed Markov (or iterated function) system}; the associated repellor (attractor, in the language of IFSs) is a Cantor set, $C$.
\item Snowflaking the Euclidean metric by the power $\alpha$, the Hausdorff dimension $s$ of $C$ becomes independent of $\alpha$.  
\item We take a skew product with covering maps on the Euclidean circle $\T$ to define a topologically cxc covering map $f: \XX_1 \to \XX_0$, $\XX_1 \subset \XX_0$ on a family of annuli; the associated repellor is $C \times \T$.
\item The map $f$ becomes a local homothety, and hence is metrically cxc.
\item A theorem of Tyson \cite[Theorem 15.10]{heinonen:analysis} implies that this metric realizes the conformal dimension, $Q$, of the cxc system $f: \XX_1 \to \XX_0$.
\eb
One motivation for this construction is that if one can find subsystems of a topologically cxc map $F: S^2 \to S^2$ conjugate to such a map $f: C \times \T \to C \times \T$, then the conformal dimension of $F$ is bounded below by that of $f$; cf. \cite{haissinsky:pilgrim:cxcIII}.

Let $\GGG$ be a directed multigraph (that is, loops of length one and multiple edges are allowed) with vertices $\{1, 2, \ldots, n\}$ and weighted edges defined as follows.  Given $(i, j) \in \{1, 2, \ldots, n\}^2$, denote by $\EEE_{ij}$ the set of edges $i \stackrel{e}{\rightarrow}j$.   For each edge $e$, suppose  $e$ is weighted by a positive integer $d(e)$.  We assume that $\GGG$ satisfies the {\em No Levy Cycle condition:}  in any cycle of edges $e_0, e_1, \ldots, d_{p-1}$ with $e_k \in \EEE_{i_ki_{k+1 \bmod p}}$, (1) $d(e_0)d(e_1)\ldots d(e_{p-1})>1$, and (2) for some $k\in\{0, \ldots p-1\}$, $\#\EEE_{i_ki_{k+1}}\geq 2$.  Furthermore, we assume that $\GGG$ is {\em irreducible}: given any pair $(i,j)$, there exists a directed path from $i$ to $j$.  

Let $\alpha>0$, and let $A_\alpha$ be the matrix given by 
\[ (A_\alpha)_{ij} = \sum_{e\in \EEE_{ij}} d(e)^{- 1/\alpha}.\]
The assumptions imply that as a function of $\alpha$, the spectral radius $\lambda(A_\alpha)$ 
is strictly monotone increasing and satisfies $\lim_{\alpha \to \infty}\lambda(A_\alpha) \geq 2$ and 
$\lim_{\alpha \to 0^+}\lambda(A_\alpha)=0$; see \cite{haissinsky:pilgrim:cxcIII} and \cite{mauldin:williams:tams88}.    

Fix $\alpha$ for which $\lambda(A_\alpha)<1$. From the theory of nonnegative matrices, it follows  that 
there exists a  
vector $w=(w_1, \ldots, w_n)$ with each $w_i>0$ such that $A_\alpha w < w$.    
For $i=1, \ldots, n$ let $I_i$ be an open  Euclidean interval of length $w_i$, and denote by $\II_0 = I_1 \sqcup \ldots \sqcup I_n$ the disjoint union of these intervals.  
Given an ordered pair $(i,j)$ and $e \in \EEE_{ij}$ let $J_e$ be an open  Euclidean interval of length $w_j d(e)^{-1/\alpha}$.  
Denote by $\II_{ij} = \sqcup_{e \in \EEE_{ij}} J_e$ the disjoint union of these intervals, and by $\II_1 = \sqcup_{i=1}^n\sqcup_{j=1}^n \II_{ij}$.  
The assumption on $w$ and the definitions imply that for each $i$, there is an embedding $\sqcup_{j=1}^n \II_{ij}\hookrightarrow I_i$ giving rise to an embedding $\II_1 \hookrightarrow \II_0$ 
satisfying the following properties: (i) it is an isometry on each interval; (ii) the closure of the image of $\II_1$ is contained in the interior of $\II_0$; (iii)  the closures of the images of distinct subintervals do not intersect.   
We fix such an embedding, and henceforth identify $\II_1$ as a subset of $\II_0$.

We extend the Euclidean metric on the intervals comprising $\II_0$ to a distance function $d(\cdot, \cdot)$ on all of $\II_0$ by setting $d(x,y)=D$ whenever $x,y$ belong to different components of $\II_0$, where $D>\frac{1}{2}\max\{w_1, \ldots, w_n\}$ is a fixed positive constant; the lower bound guarantees that the triangle inequality is satisfied.  

Define $g: \II_1 \to \II_0$ by setting, for $e \in \EEE_{ij}$, the restriction $g|_{J_e}$ to be either of the two Euclidean affine homeomorphisms sending $J_e$ onto $I_j$.   
It follows that $g|_{J_e}$ is a Euclidean similarity with ratio $d(e)^{1/\alpha} \geq 1$.   
The inverse branches of these restrictions define a so-called {\em graph-directed Markov (or iterated function) system}.  
The No Levy Cycle condition and the irreducibility condition imply that this system possesses a unique attractor (repellor, in the language of \cite[\S 2.2]{kmp:ph:cxci}), 
$C$, which is a Cantor set.  Furthermore, there exists a unique positive number $\delta \leq 1$ such that $\lambda(A_{\alpha/\delta})=1$, 
and $\delta$ coincides with the Hausdorff dimension of $C$ with respect to the metric $d$; see \cite{mauldin:williams:tams88}.  

Let $d_\alpha = d^\alpha$ be the snowflaked metric on $\II_0$.  Then the Hausdorff dimension of $C$ with respect to $d_\alpha$ is $s:=\delta/\alpha$, 
which is the unique positive parameter for which $\lambda(A_{1/s})=1$.  
Thus, while $d_\alpha$ depends on an arbitrary real parameter $\alpha$, the Hausdorff dimension of $C$ does not.  
Furthermore, for each edge $e$, the restriction $g|_{J_e}$ scales ratios of distances 
with respect to $d_\alpha$ by the factor $d(e)$, which is also independent of $\alpha$.  

Let $\T=\R/\Z$ be equipped with the Euclidean metric $d_\T$, and equip $\XX_0=\II_0 \times \T$ with the product metric $d=d_\alpha \times d_\T$, 
so that $d((x_1, t_1),(x_2, t_2)) = d_\alpha(x_1, x_2) + d_\T(t_1, t_2)$.  Define $\XX_1$ similarly.  Then $\XX_0$ is a family of right  open Euclidean cylinders, 
and  $\XX_1$ is a family of open, pairwise disjoint, essential, right subcylinders compactly contained in $\XX_0$.  
Define $f: \XX_1 \to \XX_0$ by setting, for $e \in \EEE_{ij}$, the restriction $f|_{J_e \times \T}$ to be given by 
\[ f|_{J_e \times \T}(x,t) = (g(x), d(e)t).\]
By construction, the restriction $ f|_{J_e \times \T}$ scales distances with respect to $d$ by the factor $d(e)$.  Since $f$ acts as a homothety in this metric, a covering $\mathcal{U}_0$ of $X$ by small balls in $\XX_0$ will have the property that upon setting $\mathcal{U}_n$ to be the covering obtained by components of preimages of $\mathcal{U}_0$ under $f^{-n}$, axioms (1) - (5) in the definition of metric cxc will be satisfied.  It follows that  $f: \XX_1 \to \XX_0$ is metrically coarse expanding conformal in the sense of \cite[\S 2.5]{kmp:ph:cxci} with disconnected repellor $X=C\times \T$.    By construction, the Hausdorff dimension of $X$ is $1+s$.  The  system $f: \XX_1 \to \XX_0$ again defines a unique conformal gauge and has an associated conformal dimension (technically, it is necessary to require that the conjugacies extend to the ambient spaces $\XX_0, \XX_1$; we refer to \cite[\S 2]{kmp:ph:cxci} for details).   By a theorem of Tyson \cite[Theorem 15.10]{heinonen:analysis}, the metric $d$  realizes the conformal dimension of $f$.

\section{Sierpi\'nski carpets, gaskets, and Menger spaces}

Recall that the Sierpi\'nski carpet, $S$, is the metric space obtained by starting with the unit square, subdividing into nine squares, removing the middle square, repeating with the remaining squares, and continuing; see below for an alternative description.  

On the one hand, it is well-known that there exist  hyperbolic rational maps whose Julia set is homeomorphic to $S$ (see the Appendix by Tan Lei in \cite{milnor:quadratic}).  In fact, there are many (see e.g. \cite{devaney_et_al:Sierpinski:etds}). Similarly, there exist limit sets of convex compact Kleinian groups homeomorphic to $S$.    Such examples provide a large class of metric spaces homeomorphic to $S$ and supporting a rich collection of maps which are either locally (in the case of maps) or globally quasisymmetric (in the case of groups).  

On the other hand, the $S$ is quite rigid: Bonk and Merenkov \cite{bonk:merenkov:rigidcarpet} show that the group of quasisymmetric self-maps of $S$ consists of the eight dihedral Euclidean symmetries and nothing else.  Therefore, $S$ cannot be quasisymmetrically equivalent to the boundary at infinity of any hyperbolic group.  

In contrast, Stark \cite[Theorem 2.2]{stark:menger} showed that $S$, and more generally the so-called {\em Menger spaces}, admit (in  our terminology) metrically cxc maps which, away from a thin branch locus, are homotheties with constant expansion factor.   In the remainder of this section, we briefly review Stark's construction, pose some questions, and comment on some related constructions.

\subsection*{Menger spaces} 

The following construction, and its properties, are found in \cite{stark:menger}.  
Let $n \geq 0$ be a nonnegative integer and let $k \geq 2n+1$; $I^k$ denotes the $k$-cell 
$[0,1]^k \subset \R^k$.  Let $G$ denote the subgroup of isometries of $\R^k$ generated by reflections in the faces of $I^k$, 
and $r: \R^k \to I^k$ the quotient map.  
Let $s: \R^k \to \R^k$ be given by $s(x)=3x$.  Then $s$ induces a map $f: I^k \to I^k$ on the quotient space, given by the formula $f(x)=s(r(x))$. 
Let $U_1 =\{(x_1, \ldots, x_k) \in I^k : 1/3 < x_i < 2/3 \; \mbox{for at least $n+1$ of}\; 1, 2, \ldots, k\}$.  
Set $Y_0 = I^k \setminus U_1$ and inductively put $Y_l=Y_{l-1}\setminus f^{-1}(Y_{l-1})$.  
Then $X=\cap_{l \geq 0} Y_l$ is, by a characterization theorem of Bestvina \cite[p. 2]{bestvina:menger}, the {\em Menger universal $n$-dimensional space}. 
The restriction $f: X \to X$ of $f$ to $X$ defines an open, closed, and finite-to-one map which is easily seen to be a branched covering satisfying 
axiom [Exp].  Since iterates of $s$ are all unramified, the ramification of iterates of $f: X \to X$ is uniformly bounded by the ramification of $r$, 
so axiom [Deg] is satisfied.  Given any point $x \in X$, there is a neighborhood $V$ of $x$ such that $|f(x)-f(y)| = 3|x-y|$, 
which implies immediately that the roundness and diameter distortion axioms are satisfied.  Thus $f: X \to X$ is metrically cxc.  
Note that $f$ is ramified: e.g. when $n=1, k=3$ the branch locus is the set of points in $X$ for which exactly one coordinate lies in the set $\{1/3, 2/3\}$.  

One may generalize the preceding construction.  Suppose $3 \leq \lambda_1 \leq \lambda_2 \leq \ldots \leq \lambda_n$ are integers.  Replace $s: x \mapsto 3x$ with $s: (x_1, \ldots, x_k) \mapsto (\lambda_1x_1, \ldots, \lambda_nx_n)$, put $f=r \circ s$ as before, set $\epsilon_i = \log 3/ \log \lambda_i$, and put $d(x,y) = \max_i |x_i-y_i|^{\epsilon_i}$.  Then away from the branch locus, $f: I^k \to I^k$ is again locally a homothety with factor $3$, now  with respect to the snowflaked metric $d$.  Defining $X$ as before, we get a metrically cxc dynamical system whose conformal gauge a priori may depend on the choice of expansion factors $\lambda_i$.   By replacing the middle cell $U_1$ with other suitable collections of finite cells whose closures are disjoint from each other and from the boundary of $I^k$, one obtains similar examples, with different combinatorics.  

\begin{question}
Given a degree $m \geq 2$ and a conformal gauge $\GGG$ of carpet or Menger space as above, 
how many topological conjugacy classes of metrically cxc maps $f: X \to X, X \in \GGG$, are there with $\deg(f)=m$?  
\end{question}

\subsection*{Sierpi\'nski carpets}
By taking $n=1, k=3$, and restricting to a face in a coordinate plane, one obtains a branched, metrically cxc map on $S$.

\begin{question}
Is there an unbranched metrically cxc map on the standard Sierpi\'nski carpet?  
\end{question}

If $S$ is quasisymmetrically equivalent to the Julia set of a hyperbolic rational map, then the answer is ``yes''.  If one weakens the hypothesis so as to replace $S$ with a compact metric space locally homeomorphic to $S$, then the answer is ``yes'':  
take $n=1, k=3$, and replace the group $G$ above by the group of integer translations in the coordinate directions so that the quotient space $\R^3/\Z^3$ is the three-torus $T^3$.  Then the construction above produces a set $X \subset T^3$ such that $f: 
X \to X$ is again metrically cxc with respect to the induced Euclidean metric.  The set $X$ (by \cite[Thm. 2.1]{stark:menger}) is locally homeomorphic to the {\em Menger space} of topological dimension $1$, and  the restriction of $f$ to the intersection of $X$ with the image of a coordinate plane under the natural projection yields a metric cxc system on a space $Z$ locally homeomorphic to $S$.  The space $Z$, however, cannot be embedded in plane, since the image of e.g. a coordinate axis under the natural projection will be a nonseparating simple closed curve in $Z$.

\subsection*{Sierpi\'nski gasket}

The Sierpi\'nski gasket is obtained by starting with an equilateral triangle, subdividing into four congruent equilateral triangles, removing the open middle triangle, and repeating.   Metrically cxc maps do exist on the Sierpi\'nski gasket.  Kameyama \cite{kameyama:self-similar} shows that the three Euclidean $1/3$-similitudes defining the standard triangular Sierpi\'nski gasket may be composed with rotations so that the resulting maps extend, in the manner of \S 2,  as  the inverse branches of a continuous branched covering map of the gasket to itself with three branch points; it follows easily that this yields a metrically cxc dynamical system. 
Kameyama further shows that this map, restricted to the gasket,  is topologically conjugate to the rational function $z \mapsto z^2-\frac{16}{27z}$ on its Julia set
\cite[Example 1]{kameyama:self-similar}.   By Theorem \ref{thm:promotion}, this conjugacy is quasisymmetric.    
As mentioned in \S 1, the conformal dimension of the standard Sierpi\'nski gasket is equal to one, so we conclude that the conformal dimension of the Julia set of $z \mapsto z^2-\frac{16}{27z}$ is equal to one, but not realized.  To our knowledge, this is the first nontrivial computation of a conformal dimension of a Julia set.  

This example also generalizes.  Other well-known fractals, such as the {\em hexagasket} and other {\em polygaskets}, arise as repellors of piecewise affine maps with branching as above; see e.g. \cite{eroglu:rohde:solomyak:qs}.  

\subsection*{Comment on ``flap spaces''} 

We remark that Sierpi\'nski carpets also arise naturally when the so-called ``flaps'' are excised from the ``flap spaces'' discussed by Bonk \cite{bonk:icm:qcgeom}.  However, the restriction of the dynamics to these carpets is not cxc, as this restriction fails to be an open map:  points on the boundaries of  ``holes'' have neighborhoods which map to ``half-neighborhoods''.  Thus, while of natural dynamical origin and of inherent interest (cf. \cite{merenkov:hopfian}), the associated dynamics lies outside the scope of the framework we develop in \cite{kmp:ph:cxci}.

\section{A one-parameter family of topologically cxc maps}

Let $\R^2$ denote the Euclidean plane.   
Consider the plane wallpaper group $G = \{(x,y) \mapsto \pm(x,y) + (m,n) | m, n \in \Z\} < \Isom^+(\R^2)$.  
The closure of a  fundamental domain is the rectangle $R=[0,1/2] \times [-1/2,1/2]$.  
The quotient space $\OOO=\R^2/G = R/\sim$ is homeomorphic to the two-sphere.  Away from the fixed-points of the elements of $G$, the Euclidean metric descends to a Riemannian metric on the quotient.  The completion of this metric yields a length metric $\rho$ on $\OOO$ such that the four ``corners'' (the images of the points $(0,0), (0,1/2), (1/2, 0), (1/2, 1/2)$) become cone points at which the total angle is $\pi$.
We think of $\OOO$ as a square pillowcase; in particular, it has a natural cell structure given by the vertices, edges, and faces of the two squares. 
The involution $j: \OOO \to \OOO$ induced by the map $(x,y) \mapsto (x, -y)$ on the plane descends to a map of $\OOO$ which we also denote 
suggestively by $p \mapsto \cl{p}$.   
In  $R$, this involution is reflection in the segment $\alpha$ indicated in Figure 1.

\begin{figure}[h]
\label{fig:pillowcase}
\begin{center}
\psfragscanon
\psfrag{a}{$(0,-1/2)$}
\psfrag{b}{$(1/2, -1/2)$}
\psfrag{c}{$(0,1/2)$}
\psfrag{d}{$(1/2, 1/2)$}
\psfrag{Q}{$Q$}
\psfrag{B}{$\overline{Q}$}
\psfrag{al}{$\beta$}
\psfrag{g}{$\alpha$}
\includegraphics[width=2in]{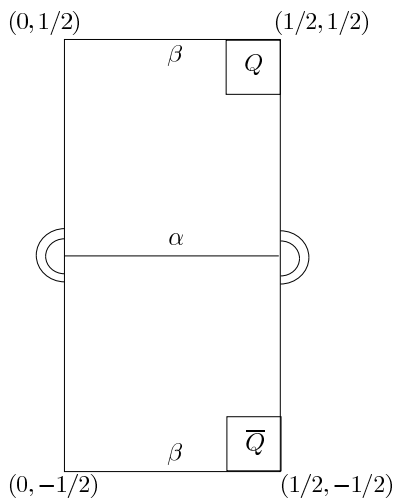}
\caption{The sphere as a square pillowcase.}  
\end{center}
\end{figure}

In particular, away from the corners, $\OOO$  inherits the coarser structure of a piecewise-affine real manifold with an  orientation-reversing symmetry.  

In this section, we show the existence of a one-parameter family $f_a: \OOO \to \OOO, a \in [0, 1/8]$ of maps 
of the sphere to itself with the following properties.  

\be
\item The map $a \mapsto f_a$ is continuous  from $[0, 1/8]$ to $C^0(\OOO, \OOO)$.
\item For all $a$, the map $f_a$ 
\be
\item is symmetric, i.e. commutes with the involution $j$;
\item  is a piecewise-affine branched covering of degree $4$ for which $f_a^{\circ 2}$ is uniformly expanding with respect to $\rho$;
\item has postcritical set given by 
\[ P_{f_a}=\{(0,0), (1/2, 0), (0, 1/2), ((1-a)/2, 1/2)\}\union \{(\tau^{\circ n}(a),0) | n \geq 0\}\]
where $\tau: [0,\frac{1}{2}] \to [0, \frac{1}{2}]$ is the full tent map given by the formula 
\[ \tau(x)=\frac{1}{2}-2|x-1/4|;\]
\item is postcritically finite if and only if $a \in \Q\intersect [0,1/8]$;  
\item is topologically cxc.
\eb
\item The map $f_0$ coincides with the integral Latt\`es map $F$ induced by $(x,y) \mapsto 2(x,y)$.
\item If $a \neq b$ then the maps $f_a, f_b$ are not topologically conjugate by any homeomorphism commuting with the involution $j$.

\item For $a \neq 0$, the map $f_a$ is not topologically conjugate to a rational function.  Let $\alpha = \{(x,0) | 0 \leq x \leq 1/2\}/\sim $ and $\beta = \{(x,1/2)| 0 \leq x \leq 1/2\}/\sim$ be the bottom and top edges, respectively, 
of $\OOO$ regarded as a square pillowcase.   
If $a>0$, the map $f_a$ has an obstruction $\Gamma = \{\gamma\}$ where 
$\gamma=\{ (x, 1/4) : 0 \leq x \leq 1/2\}/\sim  \union \{(x, -1/4) : 0 \leq x \leq 1/2\}/\sim$ is a ``horizontal'' 
simple closed curve avoiding $\alpha \union \beta$.  
\item For all $a$, there exists a homeomorphism $h_a: \OOO \to \OOO$ such that $h_a \circ f_a = f_a \circ h_a$, and $h_a$
is isotopic relative to the set $\OOO-\{\alpha \union \beta\}$ to the second iterate of a Dehn twist about the curve $\gamma$.  
\item The $1$-skeleton is forward-invariant under $f_a$; in particular, for each $a \in \Q\intersect [0,1/8]$, the map $f_a$ is the underlying map in a {\em finite subdivision rule} on the sphere in the sense of  \cite{cfp:fsr}.  
\eb

\noindent{\bf Remark:}  Since $(x,y) \mapsto 2(x,y)$ commutes with any linear map, the map $f_0$ admits many automorphisms which do not commute with $j$.  
We do not know how to rule out in general the existence of a non-symmetric conjugacy $h$ between $f_a$ and $f_b$.  
\gap

\begin{figure}
\label{fig:def_on_square}
\begin{center}
\psfragscanon
\psfrag{D1}{$\Delta_1$}
\psfrag{D2}{$\Delta_2$}
\psfrag{D3}{$\Delta_3$}
\psfrag{D1p}{$\Delta_1'$}
\psfrag{D2p}{$\Delta_2'$}
\psfrag{D3p}{$\Delta_3'$}
\psfrag{a}{$(0,0)$}
\psfrag{b}{$(a,0)$}
\psfrag{c}{$(a,a)$}
\psfrag{d}{$(0,a)$}
\includegraphics[width=3in]{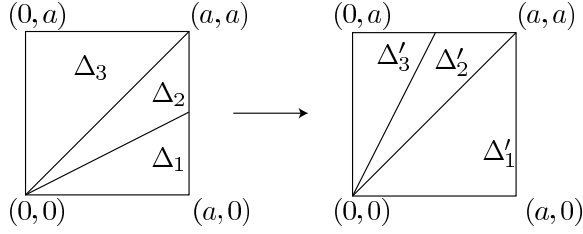}
\caption{Definition of $\tilde{R}_a$}
\end{center}
\end{figure}

\begin{figure}
\label{fig:tiling}
\begin{center}
\includegraphics[width=4in]{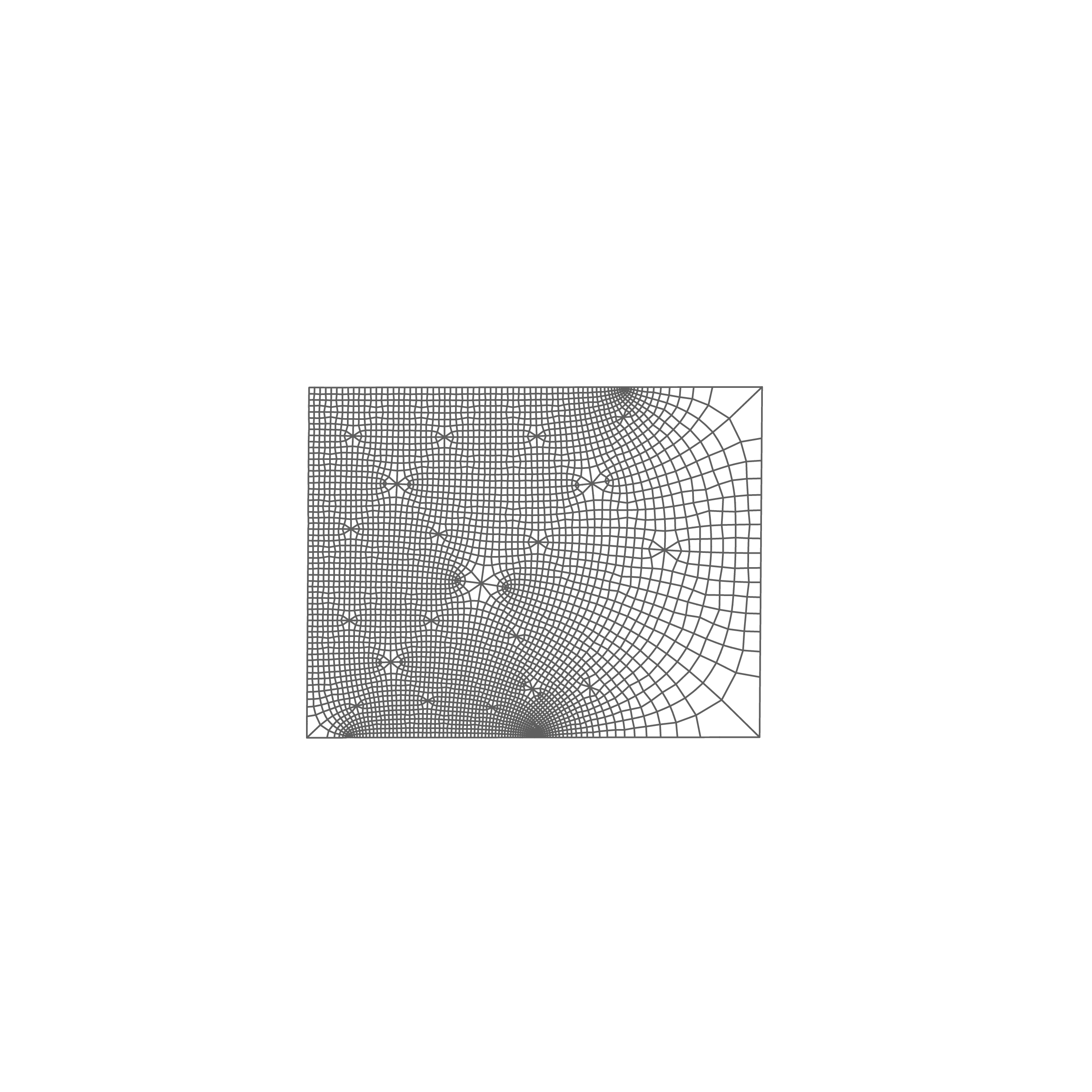}
\caption{The subdivision of the front face under the action of $f_{1/8}^{6}$ is shown.  The image is rotated by 90 degrees for convenience; circle packings are used to approximate the true conformal shape if all tiles (components of preimages of interiors of faces) are conformally equivalent to Euclidean squares (image by William Floyd).}  
\end{center}
\end{figure}

\noindent{\bf Definition of family $\mathbf{f_a}$.}  
The essential ingredient is a piecewise-linear map 
\[ \tilde{R}_a: [0,a] \times [0,a]=\Delta_1 \union \Delta_2 \union \Delta_3 \to \Delta_1' \union \Delta_2' \union \Delta_3' = [0,a] \times [0,a] \]
defined as follows.  Referring to Figure 2, 
set 
\[ \tilde{R}_a|_{\Delta_i} = T_i, \ i=1,2,3\]
where each $T_i$ is linear, and where 
\bi
\item $T_1$ is the unique linear map sending the triangle $\Delta_1$ with vertices $(0,0), (a,0), (a, a/2)$ to the triangle $\Delta_1'$ with vertices $(0,0), (a,0), (a,a)$, respectively;
\item $T_2$ is the unique linear map sending the triangle $\Delta_2$ with vertices $(0,0), (a, a/2), (a,a)$ to the triangle $\Delta_2'$ with vertices $(0,0), (a,a), (a/2, a)$;
\item $T_3$ is the unique linear map sending the triangle $\Delta_3$ with vertices $(0,0), (a,a), (0, a)$
to the triangle $\Delta_3'$ with vertices $(0,0), (a/2, a), (0,a)$.
\ib
With respect to the standard Euclidean basis, the matrices are given by
\[ T_1=\mtwo{1}{0}{0}{2}, \ T_2 = \mtwo{1}{1/2}{1}{1}\mtwo{1}{1}{1/2}{1}^{-1}, \ T_3 = \mtwo{1/2}{0}{0}{1}.\]
Recall that the singular values of a real matrix $T$ are the eigenvalues of $\sqrt{TT^t}$, 
and that the largest singular value is the Euclidean operator norm.  
The singular values of the three matrices above are given by \[ \{1,2\}, \{1/2, 2\}, \{1/2, 1\},\]
respectively.   
It follows that in the search for expansion, the worst that can happen is that some $T_i$ contracts 
the length of a tangent vector by the factor $1/2$. 

Recall that $F: \OOO \to \OOO$ is the integral Latt\`es map induced by $(x,y) \to 2(x,y)$.  
We will define $f_a = R_a \circ F$ where $R_a: \OOO \to \OOO$ is symmetric with respect to $p \mapsto \cl{p}$.  
The map $R_a$ will be the identity outside $Q_a \union \cl{Q}_a$, where 
\[ Q_a = [1/2-a, 1/2] \times [1/2-a, 1/2].\]
By symmetry, it is enough to define $R_a$ on $Q_a$.  Set 
\[ R_a = T \circ \tilde{R}_a \circ T^{-1}\]
where $T: [0,a] \times [0,a] \to Q_a$ is given by the translation $(x,y) \mapsto (x+ 1/2-a, y+ 1/2-a)$.  
This completes the definition of the family $f_a, 0 \leq a \leq 1/8$.
\gap

The remainder of this section is devoted to verification of Properties 1-6.

\noindent{\bf 1, 3, 7.}  The first two are obvious; while the latter follows immediately from the definitions given in \cite{cfp:fsr}.

\noindent{\bf 2.}  Property (a) holds by definition; we now prove (b).  
Since the smallest singular value of a $T_i$ that arises is $1/2$ and $F$ is a Riemannian homothety with expansion factor $2$, 
the differential $df_a$ does not decrease the length of tangent vectors.  Moreover, $f_a^{-1}(Q_a) \intersect Q_a=\emptyset$.  
Hence the second iterate $f^{\circ 2}_a$ expands the length of every tangent vector by a factor of at least two.   

To prove (c), note that  $f_a(\alpha \union \beta) \subset \alpha$ and 
$f_a|_{\alpha}: \alpha \to \alpha$ via $(x,0) \mapsto (\tau(x), 0)$ which is independent of $a$.  
The critical points $(1/2, \pm 1/2)$ both map under $f_a$ to $(1/2-a/2, 1/2)$ which in turn maps to $(a,0) \in \alpha$.   
For all $a$, the fate of the other critical points is the same for $f_a$ and for $F$.  
In particular, there are no recurrent or periodic critical points.  

(d) The formula for $\tau$ shows that $\tau(p/q)\equiv \pm 2p/q $ modulo $1$.  
Hence $a \in \Q$ iff the orbit of $a$ under $\tau$ is eventually periodic.  

(e) Let $\UUU_0$ be a finite covering of $\OOO$ by small Jordan domains.  
Expansion of $f_a^{\circ 2}$ implies that Axiom [Exp] holds.  
This in turn implies that the backward orbit of any point is dense in all of $\OOO$, and so $f_a$ satisfies 
Axiom [Irred].  The absence of recurrent or periodic critical points implies Axiom [Deg] is satisfied, 
so $f_a$ is topologically cxc.  

{\bf 4.}  A topological conjugacy $h$ between $f_a$ and $f_b$ which commutes with the involution $j$ 
must send the dynamically distinguished forward-invariant set $\alpha$ to itself.  
Hence $h|_\alpha$ conjugates $\tau$ to itself, $h|_\alpha=\id$, and $a=b$.

{\bf 5.}  The horizontal simple closed curve $\gamma$ has two preimages, each mapping by degree two, and each homotopic to  
$\gamma$ relative to the postcritical set of $f_a$.  If $a \neq 0$, this gives a so-called {\em Thurston obstruction} 
to the map $f_a$ being equivalent to a rational map, as was proved by Thurston for the postcritically finite case and by McMullen in general; 
see \cite{DH1}, \cite{ctm:renorm}.  

{\bf 6.}  Consider the two disjoint closed horizontal annuli $\wtA_1, \wtA_2$ where  
$\wtA_1 = [0,1/2]\times [1/8, 3/16] \union \cl{[0,1/2]\times [1/8, 3/16]}$ and 
$\wtA_2 = [0,1/2]\times [5/16, 3/8] \union \cl{[0,1/2]\times [5/16, 3/8] }$ in the fundamental domain $R$.  
Since $a\leq 1/8$ the definition of $f_a$ shows that $f_a|_{\wtA_1 \union \wtA_2}$ is independent of $a$, 
and that each annulus $\wtA_1, \wtA_2$ maps as a double cover of the annulus 
$A=[0,1/2] \times [1/4, 3/8] \union \cl{[0,1/2] \times [1/4, 3/8] }$.  

Let $h_0: \OOO \to \OOO$ be the second power of a right Dehn twist supported on $A$; 
there is a unique such $h_0$ if we require that it preserves the affine structure of $A$.  
It follows that $h_0$  lifts to a homeomorphism $h_1: \OOO \to \OOO$ which is a single right 
Dehn twist on each $\wtA_i$ and is the identity elsewhere.  It follows that $h_0$ is homotopic to $h_1$ relative 
to the top and bottom edges $\alpha \union \beta$ of the pillowcase, which contains the postcritical set $P_{f_a}$.   
Let $h_t, 0 \leq t \leq 1$, be a homotopy joining $h_0$ and $h_1$.  
By induction and homotopy lifting, we obtain a continuous family $h_t, t \geq 0$, of homeomorphisms 
such that $h_t \circ f_a = f_a \circ h_{t+1}$.  By expansion, this family is Cauchy, hence converges to a map 
$h_a: \OOO \to \OOO$ which commutes with $f_a$.  By applying the same construction to $h_0^{-1}$, we conclude that  
$h_a$ is a homeomorphism.  

\qed
\gap

\noindent{\bf Smooth versions.}  A $C^\infty$ smooth family may be constructed with similar properties as follows.  Consider a small Euclidean (cone) neighborhood $Q$ of $(1/2, 1/2)$.  Instead of $\tilde{R}_a$, use a $C^\infty$ smooth symmetric homeomorphism $\tilde{S}_a: Q \to Q$ sending $(1/2, 1/2)$ to $(1/2-a/2, 1/2)$ and which is the identity off $Q$; in suitable coordinates, one simply mollifies a small translation.  One can do this so that the differential $\tilde{S}_a$ has singular values bounded from below by a constant independent of $a$.  When $Q$ is small, the first return time to $Q$ is large.  Hence, if $a$ is sufficiently small, there is some iterate $N$ such that $f_a^{\circ N}$ is uniformly expanding.  
\gap

We do not know to what extent conjugacy classes of topologically cxc maps on manifolds contain smooth, or nearly smooth, 
representatives. 

\begin{question}
Suppose $f: S^2 \to S^2$ is topologically cxc.  Is there a smooth (smooth away from branch points, piecewise smooth, piecewise affine, \ldots ) representative in the topological conjugacy class of $f$?   
\end{question}

\subsection*{Variation of conformal dimension} 

The set of postcritically finite cxc maps $f: S^2 \to S^2$, up to topological conjugacy, is countable; hence so is the set of their corresponding conformal dimensions.  Given positive integers $m \geq n$, by looking at a map on a sphere induced by the map $(x,y) \mapsto (mx, ny)$ on the torus, and snowflaking in one direction, one can produce an example realizing the conformal dimension $1+\frac{\log m}{\log n}$.  

Beyond the postcritically finite maps, the above family shows that there exist continuous, hence uncountable, families of topologically cxc maps.  How might their conformal dimensions vary? 
In the above family, all maps $f_a, a \in (0, 1/8)$ are obstructed.   For each such map, up to homotopy $\Gamma$ is the only obstruction: a homotopically distinct obstruction $\Gamma'$ would have nontrivial geometric intersection number with $\Gamma$, and there is always an obstruction disjoint from all other obstructions \cite[Thm. 1]{kmp:canonical}.  The combinatorics of this obstruction is encoded by the matrix $(1/2 + 1/2) = (1)$, which is constant in $a$.  In \cite{haissinsky:pilgrim:cxcIII} it is shown that in general, the associated {\em snowflaked Thurston matrices} as in \S 3 give lower bounds on the Ahlfors regular conformal dimension (the corresponding infimum over all metric spaces of positive and finite Hausdorff measure in their Hausdorff dimension) of the associated metric dynamical system.  
 As a warmup, one might try to answer the following. 

\begin{question}
Does there exist a continuous, one-parameter family of obstructed cxc maps $f_t: S^2 \to S^2$, $0 \leq t \leq 1$, such that $f_0, f_1$ are (i) not topologically conjugate, and (ii) the sets of combinatorics of obstructions arising from $f_0$ and from $f_1$, as encoded by weighted directed graphs, are distinct?
 \end{question}

\def\cprime{$'$}

%

\begin{thebibliography}{BDL+}


\bibitem[Ban]{bandt:nonlinearity:mandelbrot}
Christoph Bandt.
\newblock {On the {M}andelbrot set for pairs of linear maps}.
\newblock {\em Nonlinearity} {\bf 15} (2002), 1127--1147.

\bibitem[BK]{bandt:keller:simple}
Christoph Bandt and Karsten Keller.
\newblock {Self-similar sets. {II}. {A} simple approach to the topological
  structure of fractals}.
\newblock {\em Math. Nachr.} {\bf 154} (1991), 27--39.

\bibitem[BR]{bandt:rao:topology}
Christoph Bandt and Hui Rao.
\newblock {Topology and separation of self-similar fractals in the plane}.
\newblock {\em Nonlinearity} {\bf 20} (2007), 1463--1474.

\bibitem[Be]{bestvina:menger}
Mladen Bestvina.
\newblock{Characterizing universal $k$-dimensional {M}enger compacta}.
\newblock{Memoirs of the American Mathematical Society} {\bf 380}, 1988.

\bibitem[BDL+]{devaney_et_al:Sierpinski:etds}
Paul Blanchard, Robert~L. Devaney, Daniel~M. Look, Pradipta Seal, and Yakov
  Shapiro.
\newblock {Sierpinski-curve {J}ulia sets and singular perturbations of complex
  polynomials}.
\newblock {\em Ergodic Theory Dynam. Systems} {\bf 25} (2005), 1047--1055.

\bibitem[Bon]{bonk:icm:qcgeom}
Mario Bonk.
\newblock {Quasiconformal geometry of fractals}.
\newblock In {\em International Congress of Mathematicians. Vol. II}, pages
  1349--1373. Eur. Math. Soc., Z\"urich, 2006.

\bibitem[BM]{bonk:merenkov:rigidcarpet}
Mario Bonk and Sergiy Merenkov.
\newblock {Quasisymmetric rigidity of Sierpi\'nski carpets}.
\newblock Preprint, 2009.

\bibitem[CFP]{cfp:fsr}
J.~W.~Cannon, W.~J.~Floyd, and W.~R.~Parry.
\newblock {Finite subdivision rules}.
\newblock {\em Conformal Geometry and Dynamics} (electronic), {\bf 5}(2001), 153-196.

\bibitem[DH]{DH1}
Adrien Douady and John Hubbard.
\newblock {A Proof of {T}hurston's Topological Characterization of Rational
 Functions}.
\newblock {\em Acta. Math.} {\bf 171}(1993), 263--297.


\bibitem[Edm]{edmonds:fbc}
Allan~L. Edmonds.
\newblock {Branched coverings and orbit maps}.
\newblock {\em Michigan Math. J.} {\bf 23}(1976), 289--301 (1977).

\bibitem[ERS]{eroglu:rohde:solomyak:qs}
Kemal Ero\v{g}lu, Steffen Rohde, and Boris Solomyak.
\newblock {Quasisymmetric conjugacy between quadratic dynamics and iterated
  function systems}.
\newblock To appear, J. Ergodic Th. \& Dynam. Sys., arxiv arXiv:0806.3952, v2
  2009.

\bibitem[HP1]{haissinsky:pilgrim:cxcIII}
Peter Ha{\"\i}ssinsky and Kevin Pilgrim.
\newblock {Thurston obstructions and Ahlfors regular conformal dimension}.
\newblock {\em Journal de Math\'ematiques Pures et Appliqu\'ees} {\bf 90}(2008),
  229--241.

\bibitem[HP2]{kmp:ph:cxci}
Peter Ha\"{\i}ssinsky and Kevin~M. Pilgrim.
\newblock {Coarse expanding conformal dynamics}.
\newblock  Ast\'erisque  No. 325  (2009), viii+139 pp.

\bibitem[Hei]{heinonen:analysis}
Juha Heinonen.
\newblock {\em Lectures on analysis on metric spaces}.
\newblock Universitext. Springer-Verlag, New York, 2001.

\bibitem[Kam]{kameyama:self-similar}
Atsushi Kameyama.
\newblock {Julia sets of postcritically finite rational maps and topological
  self-similar sets}.
\newblock {\em Nonlinearity} {\bf 13} (2000), 165--188.

\bibitem[MW]{mauldin:williams:tams88}
R.~Daniel Mauldin and Stanley~C. Williams.
\newblock {Hausdorff dimension in graph directed constructions}.
\newblock {\em Trans. Amer. Math. Soc.} {\bf 309} (1988), 811--829.

\bibitem[McM]{ctm:renorm}
Curtis~T. McMullen.
\newblock {\em Complex dynamics and renormalization}.
\newblock Princeton University Press, Princeton, NJ, 1994.

\bibitem[Mer]{merenkov:hopfian}
Sergiy Merenkov.
\newblock {A Sierpi\'nski carpet with the co-Hopfian property}.
\newblock {\em Invent. Math.}  {\bf 180}  (2010),  no. 2, 361--388.

\bibitem[Mil]{milnor:quadratic}
John Milnor.
\newblock {Geometry and dynamics of quadratic rational maps}.
\newblock {\em Experiment. Math.} {\bf 2} (1993), 37--83.
\newblock With an appendix by the author and Tan~Lei.

\bibitem[Pil]{kmp:canonical}
Kevin M. Pilgrim.
\newblock{Canonical Thurston obstructions}.
\newblock{\em Advances in Mathematics} {\bf 158} (2001),  no. 2, 154-168.

\bibitem[Sta]{stark:menger}
Christopher~W. Stark.
\newblock {Minimal dynamics on {M}enger manifolds}.
\newblock {\em Topology Appl.} {\bf 90} (1998), 21--30.

\bibitem[TW]{tyson:wu:selfsimilar}
Jeremy~T. Tyson and Jang-Mei Wu.
\newblock {Quasiconformal dimensions of self-similar fractals}.
\newblock {\em Rev. Mat. Iberoam.} {\bf 22} (2006), 205--258.

\end{thebibliography}
\end{document}